\documentclass[12pt]{article}
\usepackage{amssymb,amsmath}
\usepackage{graphicx}
\usepackage{url}
\usepackage{a4wide}
\usepackage{color}
\usepackage{suffix}
\newcommand {\eq} [1] {\begin{equation}\label{#1}}
\newcommand {\en} {\end{equation}}
\newcommand {\cA}       {{\mathcal A}}
\newcommand {\cB}       {{\mathcal B}}
\newcommand {\cC}       {{\mathcal C}}
\newcommand {\cD}       {{\mathcal D}}
\newcommand {\cE}       {{\mathcal E}}

\newcommand {\cL}       {{\mathcal L}}

\newcommand {\cN}       {{\mathcal N}}

\newcommand {\cQ}       {{\mathcal Q}}

\newcommand {\cS}       {{\mathcal S}}
\newcommand {\cT}       {{\mathcal T}}

\newcommand {\cW}       {{\mathcal W}}
\newcommand {\cX}       {{\mathcal X}}

\newcommand {\proof} {\par{\it Proof}. \ignorespaces}
\newcommand {\eproof}
      {\space
        {\ \vbox{\hrule\hbox{\vrule height1.3ex\hskip0.8ex\vrule}\hrule}}
        \par}
\newcommand {\C}        {{\mathbb C}}
\newcommand {\R}        {{\mathbb R}}

\newcommand {\Rn}       {\R^n}
\newcommand {\Rm}       {\R^m}
\newcommand {\Rp}       {\R^p}
\newcommand {\Rnn}      {\R^{n \times n}}
\newcommand {\Rnm}      {\R^{n \times m}}

\newcommand {\Rmm}      {\R^{m \times m}}
\newcommand {\Rmn}      {\R^{m \times n}}

\newcommand {\mat}      [1] {\left[\begin{array}{#1}}
\newcommand {\rix}          {\end{array}\right]}
\newtheorem{theorem}           {Theorem}
\newtheorem{lemma}    [theorem]{Lemma}
\newtheorem{definition}    [theorem]{Definition}
\newtheorem{corollary}[theorem]{Corollary}

\newtheorem{example}           {Example}
\newtheorem{remark}            {Remark}

\newtheorem{proposition}               {Proposition}

\newcommand {\rank}     {\mathop{\rm rank}\nolimits}

 \catcode`@=11
 \font\tenex=cmex10 
 \newdimen\p@renwd
 \setbox0=\hbox{\tenex B} \p@renwd=\wd0 
 \def\bmat#1{\begingroup \m@th
   \setbox\z@\vbox{\def\cr{\crcr\noalign{\kern2\p@\global\let\cr\endline}}%
     \ialign{$##$\hfil\kern2\p@\kern\p@renwd&\thinspace\hfil$##$\hfil
       &&\quad\hfil$##$\hfil\crcr
       \omit\strut\hfil\crcr\noalign{\kern-\baselineskip}%
       #1\crcr\omit\strut\cr}}%
   \setbox\tw@\vbox{\unvcopy\z@\global\setbox\@ne\lastbox}%
   \setbox\tw@\hbox{\unhbox\@ne\unskip\global\setbox\@ne\lastbox}%
   \setbox\tw@\hbox{$\kern\wd\@ne\kern-\p@renwd\left[\kern-\wd\@ne
     \global\setbox\@ne\vbox{\box\@ne\kern2\p@}%
     \vcenter{\kern-\ht\@ne\unvbox\z@\kern-\baselineskip}\,\right]$}%
   \null\;\vbox{\kern\ht\@ne\box\tw@}\endgroup}
 \catcode`@=12
\begin{document}

\title{Port-Hamiltonian Realizations of Positive Real Descriptor Systems}

    \author{ Chu Delin \\*
                 Department of Mathematics \\* National University of Singapore \\*
                 Singapore 119076\\*
                Email:  matchudl@nus.edu.sg.
    \and
Volker Mehrmann \\*
      Institut f\"{u}r Mathematik\\*
      MA 4-5, TU Berlin \\*
      Str. des 17. Juni 136, D-10623 Berlin, FRG \\*
      Email: mehrmann@math.tu-berlin.de.}

\maketitle
\begin{abstract} The  relationship between port-Hamiltonian and positive real descriptor systems is investigated. It is well-known that port-Hamiltonian
systems are positive real, but the converse implication does not always hold. In [K. Cherifi, H. Gernandt, and D. Hinsen. The difference between port-Hamiltonian, passive and positive real descriptor systems. Math. Control Signals Systems, pages 1–32, 2023] 
sufficient conditions for the converse are presented. We refine these conditions and present for a completely controllable, completely observable and positive real descriptor system a necessary and sufficient condition as well as an explicit method to compute  a port-Hamiltonian realization of a general positive real descriptor system.
\end{abstract}


{\bf AMS} 93B05, 93B40, 93B52, 65F35

\section{Introduction}\label{intro}
In this paper the  relationship between two  special classes of linear descriptor systems is studied. These are port-Hamiltonian  and  positive real systems. But let us first discuss general linear descriptor systems.

A linear time-invariant (LTI) \emph{descriptor system} on a time interval $\mathbb I=[t_0,t_1]$ has the form
\begin{eqnarray}
E {\dot x}(t) &=& Ax(t)+Bu(t), \nonumber  \\
         y(t) &=&  Cx(t)+Du(t)    \label{1.1}
\end{eqnarray}
with coefficient matrices $E, A\in \Rnn$, $B\in \Rnm$, $C\in \Rmn$, $D\in \Rmm$. Here  $x:\mathbb I\to \Rn$ is the state, $y:\mathbb I \to \Rp$ is the output, and  $u:\mathbb I \to \Rm$ is the input of the system. A specific property of descriptor system is that $E$ is allowed to be singular.   Although more general systems are of practical interest, see e.g. \cite{KunM06}, we require that the system (\ref{1.1}), respectively the pair of coefficients $(E,A)$, is \emph{regular}, i.e.,  that
 \begin{equation} \label{1.3}
    {\rm det}(\alpha E-\beta A)\not= 0 {\ \ \mbox {for some}}  \ \ (\alpha, \beta)\in \C^2.
\end{equation}
The \emph{generalized eigenvalues} of a regular pair $(E, A)$ are the pairs $(\alpha_j\,,\beta_j)\,\in\C^2\backslash\{(0,0)\}$ such that
\begin{equation} \label{1.3b}
    {\rm det}(\alpha_j E-\beta_j A) = 0\,, \ \ \ \ j = 1, 2,\ldots,n\,.
\end{equation}
If $\,\beta_j \not= 0$, then the eigenvalue is said to be {\em finite} with value  $\,\lambda_j = \alpha_j / \beta_j$, and otherwise, if $\,\beta_j = 0$, then the pair is  the eigenvalue {\em infinity}.
The maximal number of finite eigenvalues that a pair $(E, A)$ can have is less than or equal to the rank of $E$.

If the system (\ref{1.1}) is regular, then for \emph{consistent initial values} $ x(t_0)=x_0$ and sufficiently smooth inputs $u$, the existence and uniqueness of classical solutions is guaranteed.
The solutions are then characterized in terms of the \emph{Weierstra\ss\ Canonical Form (WCF)} of the pair $(E,A)$, which states that there exist nonsingular matrices $X$ and $Y$ such that
\eq{KCF} XEY=\bmat{ & \mu_1 & n-\mu_1 \cr
     \mu_1  &  I & 0 \cr
     n-\mu_1   & 0 & J \cr}, \quad XAY=\bmat{ & \mu_1 & n-\mu_1 \cr                              \mu_1  &  A_1 & 0 \cr
    n-\mu_1 & 0 & I \cr},
\en
where
\[ J=\mat{ccc} J_1 &           &                  \\
 & \ddots &                 \\
 &           & J_k          \rix, \ \quad
   J_i=\mat{cccc}  0 & 1        &           &      \\
    & \ddots & \ddots &       \\
    &           & 0         & 1   \\
    &           &            & 0  \rix \in \R^{n_i\times n_i}, \ i=1, \cdots, k.
\]
Here $I$ is the identity matrix of appropriate size.
The eigenvalues of  $A_1$ coincide with the finite eigenvalues of the pencil $(E, A)$ and $J$ is a nilpotent matrix
corresponding to the infinite eigenvalues. The \emph{index} of a descriptor system, denoted by ${\rm ind}(E, A)$,
is  the degree of nilpotency of the matrix $J$, i.e. the integer $i>0$, such that $J^i=0$, $J^{i-1}\not=0$.

It is well-known, see \cite{Dai89,KunM06} that the system  (\ref{1.1}) is regular and has \emph{index at most one} if and only if it has exactly $ \rank(E)$ finite eigenvalues, or equivalently that $J=0$ in (\ref{KCF}).

Assuming for simplicity that $x_0=0$ and applying a Laplace transform to (\ref{1.1}), one obtains the \emph{transfer function}
\[
\cT(s)=C(sE-A)^{-1}B+D,
\]
that maps the Laplace transform of the input to that  of the output.
The transfer function is invariant under system transformations, such as those in \eqref{KCF}, and so there is whole equivalence class of possible quintuples $(\hat E, \hat A, \hat B, \hat C, \hat D)$ with the same transfer function. This leads to the following concept of realizations.
\begin{definition}\label{def:realization}  Given a transfer function of system \eqref{1.1}, then a quintuple $(\hat E, \hat A, \hat B, \hat C, \hat D)$  with $\hat E, \hat A\in \mathbb R^{\hat n \times \hat n }$, $\hat B, \hat C^T  \in \mathbb R^{\hat n \times  m }$, $\hat D \in \mathbb R^{m \times  m }$ is called a \emph{realization} of the transfer function $\cT(s)$ or the system \eqref{1.1} if
\[
\cT(s)=C(sE-A)^{-1}B+D=\hat C (s\hat E-\hat A)^{-1}\hat B+\hat D.
\]
\end{definition}

Using the notation that for a real symmetric matrix $M$, $M\geq 0$ $(M\leq 0)$ if $M$ is positive (negative) semidefinite and introducing $\C_{+}=\{ s~| s\in \C, ~ {\rm Re}(s)>0 \}$, we have the following definition of a \emph{positive real system}, see e.g. \cite{AndV73}.
\begin{definition}\label{defPR}  System (\ref{1.1}) is called \emph{positive real} if the transfer function $\cT(s)$ is analytic in $\C_{+}$ and
\[
\cT(s)+(\cT(s))^H \geq 0 \ \mbox{\rm for all}\ s\in \C_{+}. \]
\end{definition}
%

It is not easy to check positive realness of a descriptor system, necessary and sufficient conditions for this have been derived  in \cite{ChuT08}, see also Section~\ref{S2}.

For standard state-space systems, i.e. if $E=I$, it is well-known that in the representation of a transfer function, using the Kalman decomposition, the parts that are not \emph{controllable} and \emph{observable} can be removed, leading to
a system with minimal state-space dimension, see e.g. \cite{Kai80}, called \emph{minimal realization}.
For descriptor systems the situation is more complicated, see \cite{Dai89},  and in particular \cite{BunBMN99} for the approach to construct minimal realizations in a numerically stable way.

Another class of systems that has recently received  a lot of attention are so called port-Hamiltonian (pH) systems,
\cite{SchJ14} which have been shown to provide a natural setting for control of real world physical systems and also naturally generalize to the descriptor case, see \cite{MehU23} for a recent survey presenting also a large number of applications.
In the LTI case we have the following definition.
\begin{definition} \label{def:ph}
An LTI system of the form (\ref{1.1}) is called \emph{port-Hamiltonian (pH) descriptor system} if the coefficients can be expressed as
\[
\mat{cc} A & B \\ C & D \rix =\mat{cc} (J-R)Q & G-P \\ (G+P)^TQ & S+N \rix,
\]
where $S={1\over 2}(D+D^T)$, $N={1\over 2}(D-D^T)$,
\[
\mat{cc} J & G \\ -G^T & N \rix=-\mat{cc} J & G \\ -G^T & N \rix^T,
\]
\begin{equation}\label{lagrange}
Q^TE=E^TQ\geq 0,
\end{equation}
and
\begin{equation}\label{positive} \mat{cc} Q^TRQ & Q^TP \\ P^TQ & S \rix=\mat{cc} Q^TRQ & Q^TP \\ P^TQ & S \rix^T\geq 0.
\end{equation}
The quadratic function
$\mathcal H(x)={1\over 2} x^T E^TQx$ is called the \emph{Hamiltonian} of the system.
\end{definition}
One of the major properties of pH descriptor systems is the power balance equation that replaces the conservation of energy, see e.g. \cite{MehM19,MehU23}.
\begin{theorem}\label{thm:pbe}
Consider a pH descriptor system. Then for any input $u$ the \emph{power balance equation} 
  \begin{equation}\label{eq:powerBalanceEq}
    \frac{d}{dt}\mathcal H(x(t)) = - \begin{bmatrix}
        x\\ u
    \end{bmatrix}^T\mat{cc} Q^TRQ & Q^TP \\ P^TQ & S \rix \begin{bmatrix}
        x\\ u
    \end{bmatrix}+ y^Tu
  \end{equation}
  holds along any solution $x$.
  %
%
\end{theorem}
%

Both representations of LTI descriptor systems, that in frequency domain via a positive real transfer function and that in state domain via port-Hamiltonian systems are closely related, see \cite{CheGH23} for a detailed analysis of the relationship. It is well-known, see e.g. \cite{BeaMV19,CheGH23}, that for a pH descriptor system  the transfer function is positive real, but a positive real descriptor
system does not necessarily have a port-Hamiltonian realization. To characterize when the two classes are equivalent, we discuss
the following natural questions.

\textbf{Problem 1.} Given a positive real descriptor system of the form (\ref{1.1}).  Determine necessary and sufficient conditions such that the system (\ref{1.1}) is port-Hamiltonian.

\textbf{Problem 2.} Given a positive real descriptor system of the form (\ref{1.1}).  Determine a port-Hamiltonian  realization
\begin{eqnarray}
\cE {\dot x}(t) &=& \cA x(t)+\cB u(t),  \nonumber  \\
         y(t) &=&  \cC x(t)+\cD u,    \label{1.1N}
\end{eqnarray}
of the system (\ref{1.1}) in the sense
that system (\ref{1.1N}) is port-Hamiltonian and
\[
\cT(s)=C(sE-A)^{-1}B+D=\cC(s\cE-\cA)^{-1}\cB+\cD.
\]
In  \cite{CheGH23} a partial solution to Problem 1. is  presented that relies on certain controllability and observability assumptions, see Section~\ref{S2}. In Section~\ref{S3}, we extend this partial solution and present a necessary and sufficient condition.
For Problem 2., in Section~\ref{S4}, we develop  a (numerically reliable) procedure that includes the treatment of systems that are not completely controllable or  completely observable.
Concluding remarks are given in
Section~\ref{S5}.

\section{Preliminaries}\label{S2}
In this section we present some preliminary results that are used in the remainder of the paper.
%
\begin{definition} Let $E, A \in \Rnn$, $B\in \Rnm$ and $C\in \Rmn$.

(i) The descriptor system (\ref{1.1}) (or the triplet $(E, A, B)$) is called \emph{completely controllable} if
\[ \rank \mat{cc} \alpha E-\beta A & B \rix=n,  \quad \mbox{\textrm for all} \ (\alpha, \beta)\in C^2\backslash  \{(0, 0)\}. \]

(ii) The descriptor system (\ref{1.1})  (or the triplet $(E, A, C)$) is callled \emph{completely observable} if
\[  \rank \mat{c} \alpha E-\beta A \\ C \rix=n,  \quad \mbox{\textrm for all} \ (\alpha, \beta)\in C^2\backslash  \{(0, 0)\}. \]
\end{definition}
To check whether a system is completely controllable and/or completely observable as well as other properties, numerically stable  methods have been presented in \cite{BunBMN99}, that also allow to determine a subsystem that has these properties and thus leads to a minimal representation.
For this we would need the following staircase form, see e.g. \cite{BunBMN99,ChuH99,Mim93}.
This staircase form also plays an important role for the development of a numerical procedure to compute a port-Hamiltonian realization.

\begin{lemma}\label{staircase-lemma}  Consider an LTI descriptor system of the form (\ref{1.1}). Then there exist real orthogonal matrices $U$ and $V$ such that
\begin{eqnarray}
UEV&=& \bmat{ & n_1 & n_2 & n_3 & n_4 & n_5 \cr
            n_1 & \cE & 0         & 0         & E_{14} & E_{15} \cr
            n_2 & E_{21} & E_{22} & 0         & E_{24} & E_{25} \cr
            n_3 & E_{31} & E_{32} & E_{33} & E_{34} & E_{35} \cr
            n_4 & 0         & 0         & 0          & E_{44} & E_{45} \cr
            n_5 & 0         & 0         & 0          & 0         & E_{55} \cr}, \nonumber \\
UAV&=& \bmat{ & n_1 & n_2 & n_3 & n_4 & n_5 \cr
            n_1 & \cA & 0         & 0         & A_{14} & A_{15} \cr
            n_2 & A_{21} & A_{22} & 0         & A_{24} & A_{25} \cr
            n_3 & A_{31} & A_{32} & A_{33} & A_{34} & A_{35} \cr
            n_4 & 0         & 0         & 0          & A_{44} & A_{45} \cr
            n_5 & 0         & 0         & 0          & 0         & A_{55} \cr}, \label{staircase-form} \\
UB &=& \bmat{ &     \cr
                 n_1 & \cB \cr
                 n_2 & B_2 \cr
                 n_3 & B_3 \cr
                 n_4 &  0   \cr
                 n_5 & 0    \cr}, \quad
 CV=\bmat{ & n_1 & n_2 & n_3 & n_4 & n_5 \cr
                  & \cC & 0    & 0     & C_4 & C_5 \cr}, \nonumber
\end{eqnarray}
where $E_{22}$ and $E_{44}$ are nonsingular, $sE_{33}-A_{33}$ and $sE_{55}-A_{55}$ are nonsingular for all $s\in \C$,
\[
\rank \mat{cccc} \alpha \cE-\beta \cA & 0                                       & 0             & \cB \\
\alpha E_{21}-\beta A_{21}  & \alpha E_{22}-\beta A_{22} & 0                                       & B_2 \\
\alpha E_{31}-\beta A_{31}  & \alpha E_{32}-\beta A_{32} & \alpha E_{33}-\beta A_{33} & B_3 \rix
       =n_1+n_2+n_3 \
\mbox{\textrm for all} \ (\alpha,\beta) \in \C^2\backslash \{(0,0)\},
\]
and
\[
\rank \mat{c} \alpha \cE-\beta \cA \\ \cC \rix=n_1  \
\mbox{\textrm for all} \ (\alpha,\beta) \in \C^2\backslash \{(0,0)\}.
\]
\end{lemma}

The following well-known results present conditions when an LTI descriptor system of the form (\ref{1.1}) is positive real.
\begin{theorem}\label{theorem-1}\cite{ChuT08} An LTI descriptor system of the form (\ref{1.1}) that is  completely controllable and completely observable
is positive real if and only if there exist matrices $Q\in \Rnn$ and $W\in \Rnm$ such that $Q$ and $W$ satisfy the linear matrix inequality
\begin{equation}\label{1.4}
\mat{cc} A^TQ+Q^TA & A^TW+Q^TB-C^T \\ W^TA+B^TQ-C & B^TW+W^TB-D-D^T \rix \leq 0,
\end{equation}
and, furthermore, $E^TQ\geq 0$ and $E^TW=0$.
\end{theorem}
If one drops the complete controllability and complete observability assumption then one has another characterization.
\begin{theorem}\label{theorem-2}\cite{FreJ04} An LTI descriptor system of the form (\ref{1.1}) is positive real if  there exists $Q\in \Rnn$ such that
\begin{equation}\label{inequality}
\mat{cc} A^TQ+Q^TA & Q^TB-C^T \\ B^TQ-C & -D-D^T \rix \leq 0,
\end{equation}
as well as $E^TQ\geq 0$.
\end{theorem}
A necessary condition is the following.
\begin{theorem}\label{theorem-33}\cite{FreJ04} Consider a positive real LTI descriptor system of the form (\ref{1.1}) that is completely controllable, completely observable and satisfies
\begin{equation}\label{AE}  A\, {\rm ker}(E)\subset {\rm Im}(E).
\end{equation}
If, furthermore,
\begin{equation}\label{D} D+D^T\geq  \cT(0)+(\cT(0))^T,
\end{equation}
then there exists a matrix $Q$ satisfying $E^TQ\geq 0$ and  (\ref{inequality}).
\end{theorem}


The following results (presented in a setting that is relevant for us) collect some well-known properties of LTI port-Hamiltonian descriptor systems.
\begin{proposition}\label{thm:singind}\cite{MehMW18}
Consider an LTI pH descriptor system 
for which the pair
$(E,A)=(E,(J-R)Q)$ is regular. Then the following statements hold:
\begin{enumerate}
\item[\rm (i)] If $\lambda_0\in\mathbb C$ is an eigenvalue of $(E,A)$ then $\operatorname{Re}(\lambda_0)\leq 0$.
\item[\rm (ii)] If $\omega\in\mathbb R\setminus\{0\}$ and $\lambda_0=i\omega$ is an eigenvalue of $(E,A)$, then
$\lambda_0$ is semisimple.
\item[\rm (iii)] The index of $(E,A)$ is at most two.
\end{enumerate}
\end{proposition}
Further properties are summarized in the following result, see e.g. \cite{CheGH23} for parts of these results.
\begin{theorem}\label{theorem-3} Consider an LTI descriptor system of the form (\ref{1.1}).

(a) If system (\ref{1.1}) is  pH, then $Q$ in (\ref{positive}) satisfies (\ref{inequality}).

(b) If (\ref{inequality}) holds with a nonsingular $Q$ such that $E^TQ\geq 0$, then the system (\ref{1.1}) is pH.

(c)  If the system (\ref{1.1}) is completely observable, then the matrix $Q$ in (\ref{inequality}) is nonsingular and thus the system (\ref{1.1}) is pH if and only if (\ref{inequality}) and  $E^TQ\geq 0$ hold.
\end{theorem}
\proof Parts (a) and (b) follow trivially  from simple calculations.

(c) For completeness we present a constructive proof to show that a matrix $Q$ that satisfies $E^TQ\geq 0$ and (\ref{inequality}) is nonsingular if system (\ref{1.1}) is completely observable.

Suppose that  $Q$ is singular then, using e.g. the singular value decomposition (SVD) \cite{GolV96}, there exist real orthogonal matrices $U$ and $V$ such that
\[
UQV=\bmat{ & n_1 & n_2 \cr
                n_1 &  Q_{11} & 0 \cr
                n_2 &  0          & 0 \cr},
\]
with $Q_{11}$ nonsingular. Then it follows from (\ref{inequality}) and $E^TQ\geq 0$ that
\begin{eqnarray}
UEV&=&\bmat{ & n_1 & n_2 \cr
    n_1 & E_{11} & 0 \cr
    n_2 & E_{21} & E_{22} \cr},
    \ UAV=\bmat{ & n_1 & n_2 \cr n_1  & A_{11} & 0 \cr                                    n_2   & A_{21} & A_{22} \cr}, \\
UB&=&\bmat{ &    \cr
              n_1 & B_1 \cr
              n_2 & B_2 \cr}, \
   CV=\bmat{ & n_1 & n_2 \cr
                    & C_1 & 0 \cr},
\end{eqnarray}
which contradicts the property that  system (\ref{1.1}) is completely observable. Hence, $Q$ must be nonsingular.

The converse direction follows from (a) and (b).
\eproof

In \cite{CheGH23} instead of the complete observability in  Theorem \ref{theorem-3} (c) the condition of \emph{behavioral observability}, i.e. that
$\rank \mat{c} \lambda E-A \\ C \rix=n$ for all $\lambda \in \C$, is considered. However, under this weaker condition,
the matrix $Q$ in  (\ref{inequality})  may be singular as the following example shows.
\begin{example} \label{ex:singQ}  {\rm Consider an LTI desriptor system of the form (\ref{1.1})  with
\[
E=\mat{cc} 1 & 0 \\ 0 & 0 \rix, \ A=-\mat{cc} 1 & 0 \\ 0 & 1 \rix, \
 B=\mat{c} 1 \\ 0 \rix, \ C=\mat{cc} 1 & 0 \rix, \ D=0. \]
Then $\rank \mat{c} \lambda E-A \\ C \rix=2$
for all $ \lambda \in \C$, i.e., the  system (\ref{1.1})  is behaviorally observable.  However, with
\[
Q=\mat{cc} 1 & 0 \\ 0 & 0 \rix
\]
we have
\[
\mat{cc} A^TQ+Q^TA & Q^TB-C^T \\ B^TQ-C & -D-D^T \rix=\mat{ccc} -2 & 0 & 0 \\ 0 & 0 & 0 \\ 0 & 0 & 0 \rix\leq 0,
\]
and
\[
E^TQ=\mat{cc} 1 & 0 \\ 0 & 0 \rix\geq 0.
\]
%
Moreover, it is easy to show in this case that all matrices $Q$ satisfying  (\ref{inequality}) are singular.
}
\end{example}

Theorem~\ref{theorem-33} indicates that if for a completely controllable and completely observable  LTI descriptor system (\ref{1.1}),
condition (\ref{D}) holds, then the positive realness of the system (\ref{1.1}) implies that it is a pH descriptor system. However,  condition (\ref{D}) is not satisfied in general.

\section{Solution of Problem 1.}\label{S3}
In this section we study the solution of Problem 1. For this, we need two condensed forms.
\begin{lemma}\label{lemma-1a} Consider an LTI descriptor system of the form (\ref{1.1}) that is completely controllable, completely observable, and positive real.
Then there exist orthogonal matrices $U, V\in \Rnn$, and $W\in \Rmm$ such that
\begin{eqnarray}\label{EA1}
U E V &=& \bmat{ & n_1 & n_2 & n_2 & n_3 \cr
n_1  &  E_{11}     & E_{12}     & E_{13}    & E_{14}   \cr
n_2  &  0    & 0     & E_{23}    & 0    \cr
n_2  &  0    & 0     & 0    & 0 \cr
n_3 &  0 & 0 & 0 & 0 \cr},    \quad
U A V =  \bmat{ & n_1 & n_2 & n_2 & n_3 \cr
n_1 & A_{11}     & A_{12}     & A_{13}   & A_{14}     \cr
n_2 & 0    & A_{22}     & A_{23}      & A_{24} \cr
n_2 & 0    & 0     & A_{33}     &  0  \cr
n_3 & 0 & 0 & A_{43} & A_{44} \cr},\\
\label{B1}
   UBW &=& \bmat{ & m_1 & n_2 & n_3  \cr
        n_1 &  \hat B_{11} & \hat B_{12}  & \hat B_{13} \cr
        n_2 &  \hat B_{21} & \hat B_{22}  & \hat B_{23}   \cr
        n_2 &  0            & B_{32}      & B_{33} \cr
        n_3 &  0 & 0 & B_{43} \cr},
\end{eqnarray}
where $E_{11}$, $E_{23}$, $A_{22}$, $A_{33}$, $A_{44}$, $B_{32}$ and $B_{43}$ are nonsingular.
\end{lemma}
\proof We construct (\ref{EA1})
by the following numerically implementable procedure.

Step 1: Since the pencil $(E, A)$ is regular, we can compute the generalized upper triangular form of the pencil $(E, A)$, see  \cite{DemK93a,DemK93b}, that determines orthogonal matrices $U_1, V_1\in \Rnn$  such that
\[
U_1 (\lambda E-A) V_1=\bmat{ & n_1 & n-n_1 \cr
n_1  & \lambda E_{11}-A_{11} & \lambda E_{12}^{(1)}-A_{12}^{(1)} \cr
n-n_1  & 0                                  & \lambda E_{22}^{(1)}-A_{22}^{(1)} \cr},
\]
where $E_{11}$ is nonsingular, and
$ \rank (\lambda E_{22}^{(1)}-A_{22}^{(1)})=n-n_1$ for all $ \lambda \in \C$.

Step 2: Compute, for example via the SVD, orthogonal matrices $U_2, V_2 \in \R^{(n-n_1)\times (n-n_1)}$ such that
\[
U_2 E_{22}^{(1)} V_2 =\bmat{ & n_2 & n-n_1-n_2 \cr
n_2   & E_{22}^{(2)} & 0 \cr
n-n_1-n_2 &  0      & 0 \cr},
\]
where $n_2=\rank(E_{22}^{(1)})$ and $E_{22}^{(2)}$ is nonsingular. Set
\[
U_2 A_{22}^{(1)} V_2=\bmat{ & n_2 & n-n_1-n_2 \cr
n_2     & A_{22}^{(2)} & A_{23}^{(2)} \cr
n-n_1-n_2 & A_{32}^{(2)} & A_{33}^{(2)} \cr}.
\]

Step 3: Compute, for example via the SVD, orthogonal matrices $U_3, V_3\in \R^{(n-n_1-n_2)\times (n-n_1-n_2)}$ such that
\[
U_3 A_{33}^{(2)} V_3 =\bmat{ & \nu_2 & n_3 \cr
\nu_2   & 0 & 0 \cr
n_3     & 0 & A_{44} \cr},
\]
where $\nu_2=n-n_1-n_2-n_3$ and $A_{44}$ is nonsingular. Set
\begin{eqnarray*}
A_{23}^{(2)} V_3&=&\bmat{ & \nu_2 & n_3 \cr
& A_{23}^{(3)} & A_{24}^{(3)} \cr}, \quad
U_3    A_{32}^{(2)} =\bmat{ &        \cr
\nu_2 &    A_{32}^{(3)} \cr
n_3  & A_{42}^{(3)}   \cr}, \\
 E_{12}^{(1)}V_2 \mat{cc} I & \\ & V_3 \rix    &=&\bmat{ & n_2 & \nu_2 & n_3 \cr         & E_{12}^{(3)} & E_{13}^{(3)} & E_{14}^{(3)} \cr}, \\
A_{12}^{(1)}V_2 \mat{cc} I & \\ & V_3 \rix &=&\bmat{ & n_2 & \nu_2 & n_3 \cr            & A_{12}^{(3)} & A_{13}^{(3)} & A_{14}^{(3)} \cr}.
\end{eqnarray*}
Then we have
\begin{eqnarray*}
 \mat{cc} I & \\ & U_3 \rix \mat{cc} I & \\ & U_2 \rix U_1 E V_1 \mat{cc} I & \\ & V_2 \rix \mat{cc} I & \\ & V_3 \rix
   &=&\bmat{ & n_1 & n_2 & \nu_2 & n_3 \cr
         n_1  & E_{11} & E_{12}^{(3)} & E_{13}^{(3)} & E_{14}^{(3)} \cr
         n_2  & 0          & E_{22}^{(2)} & 0                  & 0                   \cr
       \nu_2 & 0          & 0                   & 0                  & 0                   \cr
         n_3  & 0          & 0                   & 0                  & 0                   \cr}, \\
         \mat{cc} I & \\ & U_3 \rix \mat{cc} I & \\ & U_2 \rix U_1 A V_1 \mat{cc} I & \\ & V_2 \rix \mat{cc} I & \\ & V_3 \rix
   &=&\bmat{ & n_1 & n_2 & \nu_2 & n_3 \cr
         n_1  & A_{11} & A_{12}^{(3)} & A_{13}^{(3)} & A_{14}^{(3)} \cr
         n_2 & 0 & A_{22}^{(2)} & A_{23}^{(3)}  & A_{24}^{(3)} \cr
      \nu_2 & 0 & A_{32}^{(3)} & 0 & 0 \cr
         n_2 & 0 & A_{42}^{(3)} & 0 & A_{44} \cr}.
\end{eqnarray*}
System  (\ref{1.1}) is completely controllable, completely observable, and positive real, so the index of the pencil $(E, A)$ is at most two, see e.g. \cite{ChuT08}. Consequently, the $J_i$ blocks in the \emph{Weierstra\ss\ Canonical Form (WCF)} of $(E, A)$ have size at most $2$ and the number of $2\times 2$ blocks is $n_2$.  Obviously, the regularity of the pencil $(E, A)$ gives that
$A_{23}^{(3)}$ is of full column rank and $A_{32}^{(3)}$ is of full row rank and so, $\nu_2\leq n_2$. Hence, we must have
$n_2=\nu_2$
and thus, $A_{23}^{(3)}$ and $A_{32}^{(3)}$ are nonsingular.

Step 4: Set
\begin{eqnarray*}
U&=& \mat{cc} I & \\ & U_3 \rix \mat{cc} I & \\ & U_2 \rix U_1, \\  V&=&V_1 \mat{cc} I & \\ & V_2 \rix \mat{cc} I & \\ & V_3 \rix \mat{cccc} I_{n_1} &&&   \\
& & I_{n_2} & \\& I_{n_2} &  & \\            & & & I_{n_3}               \rix,\\
 E_{12}&=&E_{13}^{(3)}, \ E_{13}=E_{12}^{(3)}, \ E_{14}=E_{14}^{(3)}, \ E_{23}=E_{22}^{(2)}, \\ A_{12}&=&A_{13}^{(3)}, \ A_{13}=A_{12}^{(3)}, \ A_{14}=A_{14}^{(3)}, \ A_{23}=A_{22}^{(2)}, \
     A_{22}=A_{23}^{(3)}, \ A_{24}=A_{24}^{(3)}, \\ A_{33}&=&
     A_{32}^{(3)}, \ A_{43}=A_{42}^{(3)}.
     \end{eqnarray*}
Then (\ref{EA1}) holds.

Step 5: Set
\[
UB=\bmat{ &    \cr
n_1 & B_1^{(3)} \cr
n_2 & B_2^{(3)} \cr
n_2 & B_3^{(3)} \cr
n_3 & B_4^{(3)} \cr}.
\]
Since $\mat{cc} E & B \rix$ is of full row rank, we have that $ \mat{c} B_3^{(3)} \\ B_4^{(3)} \rix$ is of full row rank.
Hence, using e.g. a permuted $QR$ decomposition, we can compute an orthogonal matrix $W$ such that
\[
\mat{c} B_3^{(3)} \\ B_4^{(3)} \rix W=\bmat{ & m_1 & n_2 & n_3 \cr
n_2      & 0  & B_{32} & B_{33} \cr
n_3       & 0  & 0         & B_{43} \cr},
\]
where $B_{32}$ and $B_{43}$ are nonsingular. Setting
\[
\mat{c} B_1^{(3)} \\ B_2^{(3)} \rix W=\bmat{ & m_1 & n_2 & n_3 \cr                        n_1 & \hat B_{11} & \hat B_{12} & \hat B_{13} \cr                                   n_2 & \hat B_{21} & \hat B_{22} & \hat B_{23} \cr}.
\]
we have obtained (\ref{B1}).
\eproof
The forms (\ref{EA1}) and (\ref{B1}) are computed by orthogonal transformations. Indices $n_1$ and $n_2$ in (\ref{EA1}) are the dimension of the generalized finite eigenspace of $(E, A)$ and the number of $2\times 2$ blocks in the  WCF of $(E, A)$, respectively. Hence,
the forms (\ref{EA1}) and (\ref{B1}) including indices $n_1$ and $n_2$ can be computed in a numerically stable way.

The following condensed form is a refined version of the forms (\ref{EA1}) and (\ref{B1}) by using nonsingular but non-orthogonal transformations.
\begin{lemma}\label{lemma-1} Consider an LTI  descriptor system of the form (\ref{1.1}) that is completely controllable, completely observable, and positive real.
Then there exist nonsingular matrices $X, Y\in \Rnn$, and $Z\in \Rmm$ such that
\begin{eqnarray}
X E Y &=&  \bmat{ & n_1 & n-n_1 \cr
                   n_1    & I & 0 \cr
               n-n_1     & 0 & \cN \cr}
       =\bmat{ & n_1 & n_2 & n_2 & n_3 \cr
             n_1  &  I     & 0     & 0    & 0   \cr
             n_2  &  0    & 0     & I    & 0    \cr
             n_2  &  0    & 0     & 0    & 0 \cr
             n_3 &  0 & 0 & 0 & 0 \cr},    \nonumber \\
X A Y &=& \bmat{ & n_1 & n-n_1 \cr
                n_1      & A_1 & 0 \cr
             n-n_1      & 0 & I \cr}
=\bmat{ & n_1 & n_2 & n_2 & n_3 \cr
       n_1 & A_1     & 0     & 0   & 0     \cr
       n_2 & 0    & I     & 0      & 0 \cr
       n_2 & 0    & 0     & I     &  0  \cr
       n_3 & 0 & 0 & 0 & I \cr}, \nonumber \\
X B Z &=& \bmat{ &     \cr
                   n_1   & B_1  \cr
              n-n_1     &  \cB \cr}
= \bmat{ & m_1 & n_2 & n_3  \cr
        n_1 &  B_{11} & B_{12}  & B_{13} \cr
        n_2 &  B_{21} & 0        & 0   \cr
        n_2 &  0            & I      & 0 \cr
        n_3 &  0 & 0 & I \cr}, \label{form-1} \\
Z^TCY &=&\bmat{ & n_1 & n-n_1 \cr
                            & C_1 & \cC \cr}
= \bmat{ & n_1 & n_2 & n_2 & n_3 \cr
        m_1 & C_{11} & C_{12} & C_{13}  & C_{14} \cr
        n_2  & C_{21} & C_{22} & C_{23}  & C_{24} \cr
        n_3  &  C_{31} & C_{32} & C_{33} & C_{34} \cr}. \nonumber
\end{eqnarray}
Moreover,
\begin{equation}\label{N}  \cC\cN\cB=\mat{ccc} 0 & C_{12} & 0 \\ 0 & C_{22} & 0 \\  0 & C_{32} & 0 \rix\leq 0, \ i.e., \ C_{12}=0, \ C_{32}=0, \ C_{22}\leq 0,
\end{equation}
and there exists a matrix $\cQ_{11}\in \R^{n_1\times n_1}$ with  $\cQ_{11}\geq 0$ such that
\begin{equation}\label{lmi-1}
\mat{cccc} \cQ_{11}^TA_1+A_1^T\cQ_{11} & \cQ_{11}^TB_{11}-C_{11}^T & \cQ_{11}^TB_{12}-C_{21}^T & \cQ_{11}^T B_{13}-C_{31}^T\\
            B_{11}^T\cQ_{11}-C_{11}         &   -S_{11}                 &         -S_{12}+C_{13}+B_{21}^TC_{22}^T & -S_{13}+C_{14}   \\
            B_{12}^T\cQ_{11}-C_{21}         &  -S_{12}^T+C_{13}^T+C_{22}B_{21}   &         -S_{22}+C_{23}+C_{23}^T & -S_{23}+C_{24}+C_{33}^T \\
            B_{13}^T\cQ_{11}-C_{31}        & -S_{13}^T+C_{14}^T & -S_{23}^T+C_{33}+C_{24}^T & -S_{33}+C_{34}+C_{34}^T \rix \leq 0,
\end{equation}
where
\[
\mat{ccc} S_{11} & S_{12} & S_{13} \\ S_{12}^T & S_{22} & S_{23} \\ S_{13}^T & S_{23}^T & S_{33}  \rix=Z^T(D+D^T)Z.
\]
\end{lemma}
\proof   We can obtain nonsingular matrices $X_1$, $Y_1$ and $Z_1$ by applying block Gaussian elimination to (\ref{EA1}) and (\ref{B1})
using the nonsingular matrices $E_{11}$, $A_{22}$, $A_{33}$, $A_{44}$, $B_{32}$ and $B_{43}$
such that
\[
X_1U (sE-A)V Y_1 =\bmat{ & n_1 & n_2 & n_2 &  n_3 \cr
    n_1   & sI-A_1 & 0 & 0 & 0 \cr
    n_2   & 0 & -I & sI & 0 \cr
    n_2   & 0 & 0 & -I & 0 \cr
    n_3   & 0 & 0 & 0 & -I \cr},
\]
and
\[
X_1 U B W Z_1= \bmat{ & m_1 & n_2 & n_3 \cr
 n_1  &  B_{11} & B_{12} & B_{13} \cr
n_2  &  B_{21} & B_{22} & B_{23} \cr
n_2 &   0         & I           & 0         \cr
n_3 &   0         & 0           & I        \cr}.
\]
Then, using block-Gaussian elimination again, we can compute nonsingular matrices $X_2$ and $Y_2$ such that
\begin{eqnarray*}
X_2    \mat{cccc} sI-A_1 & 0 & 0 & 0 \\
        0 & -I & sI & 0 \\
        0 & 0 & -I & 0 \\
    0 & 0 & 0 & -I \rix Y_2 &=&\mat{cccc} sI-A_1 & 0 & 0 & 0 \\
    0 & -I & sI & 0 \\
    0 & 0 & -I & 0 \\
    0 & 0 & 0 & -I \rix,
\\
X_2  \mat{ccc} B_{11} & B_{12} & B_{13} \\ B_{21} & B_{22} & B_{23} \\ 0  & I & 0 \\ 0 & 0 & I  \rix
  &=& \mat{ccc} B_{11} & B_{12} & B_{13} \\ B_{21} & 0 & 0 \\ 0 & I & 0 \\ 0 & 0 & I  \rix.
\end{eqnarray*}
Setting $X=X_2X_1U$, $Y=VY_1Y_2$ and $Z=WZ_1$ yields the condensed form (\ref{form-1}).
Moreover, we have
\begin{eqnarray*}
 Z^T\cT(s)Z &=& Z^T(C(sE-A)^{-1}B+D)Z  \\
  &=& \mat{c} C_{11} \\ C_{21} \\ C_{31} \rix (sI-A_1)^{-1} \mat{ccc} B_{11} & B_{12} & B_{13} \rix
        +  (Z^TDZ-\cC\cB)-s\cC\cN\cB.
\end{eqnarray*}
Since the system is still positive real, it follows that
\[
\mat{c} C_{11} \\ C_{21} \\ C_{31} \rix (sI-A_1)^{-1} \mat{ccc} B_{11} & B_{12} & B_{13} \rix
        +  (Z^TDZ-\cC\cB)
\]
is positive real and $-\cC\cN\cB\geq 0$. Hence (\ref{N}) and (\ref{lmi-1}) hold. \eproof
%


\begin{remark} \label{rem:strong}{\rm
Condition (\ref{D}) holds if and only if
$D+D^T\geq D+D^T-\cC\cB-(\cC\cB)^T$,
or equivalently,
$\cC\cB+(\cC\cB)^T\geq 0$,
which means that
\[
C_{14}=0, \ C_{13}+B_{21}^TC_{22}^T=0, \
   \mat{cc} C_{23}+C_{23}^T & C_{24}+C_{33}^T \\ C_{24}^T+C_{33} & C_{34}+C_{34}^T \rix \geq 0.
   \]
This shows that condition (\ref{D}) is very strong and will not hold in general.
}
\end{remark}

After these preparations we are now ready to present a necessary and sufficient condition that solves Problem~1.

\begin{theorem}\label{theorem-4} Consider an LTI descriptor system of the form (\ref{1.1}) that is completely controllable, completely observable, and positive real. Then it is port-Hamiltonian if and only if $D+D^T\geq 0$.
\end{theorem}
\proof  It follows from Theorem \ref{theorem-3} that  (\ref{inequality}) holds, since the system (\ref{1.1}) is port-Hamiltonian, then $D+D^T\geq 0$. Hence, the necessity follows.

To prove the sufficiency, observe that since the  system (\ref{1.1}) is completely controllable, completely observable, and positive real, by Lemma \ref{lemma-1},  the matrix inequality (\ref{lmi-1}) holds with $\cQ_{11}\geq 0$.
Additionally,
\[
\mat{ccc} -S_{11} & -S_{12}  & -S_{13} \\ -S_{12}^T & -S_{22} & -S_{23} \\ -S_{13}^T & -S_{23}^T & -S_{33} \rix=-Z^T(D+D^T)Z\leq 0,
\]
or equivalently
\[   \mat{ccc} -S_{22} & -S_{23} & -S_{12}^T \\ -S_{23}^T & -S_{33} & -S_{13}^T \\ -S_{12} & -S_{13} & -S_{11} \rix \leq 0.
\]
Hence there exist matrices $K_1$ and $K_2$ such that
\[
\mat{cc} S_{12} & S_{13} \rix =S_{11} \mat{cc} K_1 & K_2 \rix, \]
and
\[ \mat{cc} -\cS_{22} & -\cS_{23} \\
                 -\cS_{23}^T & -\cS_{33} \rix:=
 \mat{cc} -S_{22} & -S_{23} \\ -S_{23}^T & -S_{33} \rix +\mat{c} S_{12}^T \\ S_{13}^T \rix \mat{cc} K_1 & K_2 \rix \leq 0.
 \]
Choosing
\begin{eqnarray*}
\mat{cc} Q_{33} & Q_{34} \\ Q_{43} & Q_{44} \rix&=&
\mat{cc} -S_{22}+C_{23} & -S_{23}+C_{24} \\ -S_{23}^T+C_{33} & -S_{33}+C_{34} \rix\\
&&-\mat{cc} K_1 & K_2 \rix^T\mat{cc} -S_{12}+B_{21}^TC_{22}+C_{13} & -S_{13}+C_{14} \rix,
\end{eqnarray*}
and using (\ref{lmi-1}), it follows that
{\footnotesize
\begin{eqnarray*}
&&\mat{cc|ccc}
-S_{22}                 & -S_{23} & -S_{12}^T & -S_{22}+C_{23}-Q_{33} & -S_{23}+C_{24}-Q_{34} \\
-S_{23}^T             & -S_{33} & -S_{13}^T & -S_{23}^T+C_{33}-Q_{43} & -S_{33}+C_{34}-Q_{44}\\  \hline
-S_{12} & -S_{13} & -S_{11} & -S_{12}+B_{21}^TC_{22}^T+C_{13} & -S_{13}+C_{14}  \\
-S_{22}+C_{23}^T-Q_{33}^T     & -S_{23}+C_{33}^T-Q_{43}^T &  -S_{12}^T+C_{22}B_{21}+C_{13}^T & -S_{22}+C_{23}+C_{23}^T& -S_{23}+C_{24}+C_{33}^T \\
-S_{23}^T+C_{24}^T-Q_{34}^T & -S_{33}+C_{34}^T-Q_{44}^T & -S_{13}^T+C_{14}^T & -S_{23}^T+C_{24}^T+C_{33} & -S_{33}+C_{34}+C_{34}^T \rix \\
&& =  \mat{ccccc} I  &   &  &   &  \\           & I        &       &    &      \\
-K_1 & -K_2 & I    &           &   \\
&            &     & I         &     \\
&             &     &    & I \rix^{-T} \\
& & \times \mat{cc|ccc}
 -\cS_{22}  & -\cS_{23} & 0 & 0 & 0\\
-\cS_{23}^T   & -\cS_{33} & 0  & 0 & 0  \\  \hline
 0 & 0 &  -S_{11}&  -S_{12}+B_{21}^TC_{22}^T+C_{13} & -S_{13}+C_{14}  \\
0  & 0 &  -S_{12}^T+C_{22}B_{21}+C_{13}^T & -S_{22}+C_{23}+C_{23}^T& -S_{23}+C_{24}+C_{33}^T \\
 0  & 0 & -S_{13}^T+C_{14}^T & -S_{23}^T+C_{24}^T+C_{33} & -S_{33}+C_{34}+C_{34}^T \rix   \\
& &  \times\mat{ccccc} I & & & &      \\
& I        &       &          &      \\
-K_1 & -K_2 & I    &           &   \\
&            &     & I         &     \\
&             &     &           & I \rix^{-1}  \leq 0,
\end{eqnarray*}
}
or equivalently,
{\footnotesize
\begin{equation}
    \mat{ccccc}
    Q_{33}+Q_{33}^T & Q_{34}+Q_{43}^T & -C_{22}B_{21}-C_{13}^T & -Q_{33}-C_{23}^T & -Q_{34}-C_{33}^T\\
    Q_{43}+Q_{34}^T & Q_{44}+Q_{44}^T & -C_{14}^T & -Q_{43}-C_{24}^T & -Q_{44}-C_{34}^T \\
    -B_{21}^TC_{22}^T-C_{13} & -C_{14} & -S_{11} & -S_{12}+B_{21}^TC_{22}^T+C_{13} & -S_{13}+C_{14}^T  \\
    -Q_{33}^T-C_{23} & -Q_{43}^T-C_{24} &  -S_{12}^T+C_{22}B_{21}+C_{13}^T & -S_{22}+C_{23}+C_{23}^T& -S_{23}+C_{24}+C_{33}^T \\
    -Q_{34}^T-C_{33} & -Q_{44}^T-C_{34} & -S_{13}^T+C_{14}^T & -S_{23}^T+C_{24}^T+C_{33} & -S_{33}+C_{34}+C_{34}^T \rix \leq 0. \label{LMI-3}
\end{equation}
}
Hence, there exists a matrix with partitioning $\mat{cc} \cX_{11} & \cX_{12} \\ \cX_{21} & \cX_{22} \\ \cX_{31} & \cX_{32} \rix$
such that
\begin{eqnarray*} && \mat{cc} -B_{21}^TC_{22}^T-C_{13} & -C_{14}   \\
 -Q_{33}^T-C_{23} & -Q_{43}^T-C_{24}  \\
 -Q_{34}^T-C_{33} & -Q_{44}^T-C_{34}\rix \\
&&    =\mat{ccc} -S_{11} & -S_{12}+B_{21}^TC_{22}^T+C_{13} & -S_{13}+C_{14}^T  \\
-S_{12}^T+C_{22}B_{21}+C_{13}^T & -S_{22}+C_{23}+C_{23}^T& -S_{23}+C_{24}+C_{33}^T \\
-S_{13}^T+C_{14}^T & -S_{23}^T+C_{24}^T+C_{33} & -S_{33}+C_{34}+C_{34}^T \rix \mat{cc} \cX_{11} & \cX_{12} \\ \cX_{21} & \cX_{22} \\ \cX_{31} & \cX_{32} \rix,
\end{eqnarray*}
and
\begin{eqnarray*}
\Lambda
&:=& \mat{cc}  Q_{33}+Q_{33}^T & Q_{34}+Q_{43}^T \\
                    Q_{43}+Q_{34}^T & Q_{44}+Q_{44}^T \rix\\
&-& \mat{ccc} -C_{22}B_{21}-C_{13}^T & -Q_{33}-C_{23}^T & -Q_{34}-C_{33}^T\\
                   -C_{14}^T & -Q_{43}-C_{24}^T & -Q_{44}-C_{34}^T \rix \mat{cc} \cX_{11} & \cX_{12} \\ \cX_{21} & \cX_{22} \\ \cX_{31} & \cX_{32} \rix  \leq 0.
\end{eqnarray*}
Let
$Q_{11}=\cQ_{11}$ with $\cQ_{11}\geq 0$ as in (\ref{lmi-1}) and
\[
\mat{c} Q_{31} \\ Q_{41} \rix
  =\mat{cc} \cX_{11} & \cX_{12} \\ \cX_{21} & \cX_{22} \\ \cX_{31} & \cX_{32} \rix^T \mat{c} B_{11}^TQ_{11}-C_{11}   \\
B_{12}^TQ_{11}-C_{21}   \\
B_{13}^TQ_{11}-C_{31} \rix.
\]
Then it follows from  (\ref{lmi-1}) and $ \Lambda\leq 0$ that
{\scriptsize
\begin{eqnarray*}
& &  \mat{cccccc} I      &  &  & &  &   \\
&   I        &&   &&  \\
& & I  & & &  \\
& -\cX_{11} & -\cX_{12} & I &&   \\
& -\cX_{21} & -\cX_{22} &  & I & \\
& -\cX_{31} & -\cX_{32} &                  && I   \rix^T  \\
&\times&
\mat{cccccc} A_1^TQ_{11}+Q_{11}^TA_1  & Q_{31}^T & Q_{41}^T & Q_{11}^TB_{11}-C_{11}^T &   Q_{11}^TB_{12}-C_{21}^T & Q_{11}^TB_{13}-C_{31}^T\\
Q_{31}  & Q_{33}+Q_{33}^T & Q_{34}+Q_{43}^T & -C_{22}B_{21}-C_{13}^T & -Q_{33}-C_{23}^T & -Q_{34}-C_{33}^T\\
Q_{41} & Q_{43}+Q_{34}^T & Q_{44}+Q_{44}^T & -C_{14}^T & -Q_{43}-C_{24}^T & -Q_{44}-C_{34}^T \\
B_{11}^TQ_{11}-C_{11}  & -B_{21}^TC_{22}^T-C_{13} & -C_{14} & -S_{11} & -S_{12}+B_{21}^TC_{22}^T+C_{13} & -S_{13}+C_{14}^T  \\
B_{12}^TQ_{11}-C_{21}  & -Q_{33}^T-C_{23} & -Q_{43}^T-C_{24} &  -S_{12}^T+C_{22}B_{21}+C_{13}^T & -S_{22}+C_{23}+C_{23}^T& -S_{23}+C_{24}+C_{33}^T \\
B_{13}^TQ_{11}-C_{31} & -Q_{34}^T-C_{33} & -Q_{44}^T-C_{34} & -S_{13}^T+C_{14}^T & -S_{23}^T+C_{24}^T+C_{33} & -S_{33}+C_{34}+C_{34}^T \rix     \\
&\times &  \mat{cccccc} I      &            &    & &                     & \\
&   I        &  &  &  &  \\
&            & I &    &     &   \\
 & -\cX_{11} & -\cX_{12} & I                &   & \\
& -\cX_{21} & -\cX_{22} &                 & I &  \\
& -\cX_{31} & -\cX_{32} &                  &   & I  \rix  \\
&=&
\mat{ccccc} A_1^T\cQ_{11}+\cQ_{11}^TA_1  & 0              & \cQ_{11}^TB_{11}-C_{11}^T &   \cQ_{11}^TB_{12}-C_{21}^T & \cQ_{11}^TB_{13}-C_{31}^T\\
0  & \Lambda   &  0 & 0 & 0 \\
B_{11}^T\cQ_{11}-C_{11}  &  0 & -S_{11} & -S_{12}+B_{21}^TC_{22}^T+C_{13} & -S_{13}+C_{14}^T  \\
B_{12}^T\cQ_{11}-C_{21}  &  0 &   -S_{12}^T+C_{22}B_{21}+C_{13}^T & -S_{22}+C_{23}+C_{23}^T& -S_{23}+C_{24}+C_{33}^T \\
B_{13}^T\cQ_{11}-C_{31} & 0  & -S_{13}^T+C_{14}^T & -S_{23}^T+C_{24}^T+C_{33} & -S_{33}+C_{34}+C_{34}^T \rix\leq 0,
\end{eqnarray*}
}
which implies that
{{\scriptsize
\begin{equation}
\mat{cccccc} A_1^TQ_{11}+Q_{11}^TA_1  & Q_{31}^T & Q_{41}^T & Q_{11}^TB_{11}-C_{11}^T & Q_{11}^TB_{12}+Q_{31}^T-C_{21}^T & Q_{11}^TB_{13}+Q_{41}^T-C_{31}^T\\
Q_{31}  & Q_{33}+Q_{33}^T & Q_{34}+Q_{43}^T & -C_{22}B_{21}-C_{13}^T & Q_{33}^T-C_{23}^T & Q_{43}^T-C_{33}^T\\
Q_{41} & Q_{43}+Q_{34}^T & Q_{44}+Q_{44}^T & -C_{14}^T & Q_{34}^T-C_{24}^T & Q_{44}^T-C_{34}^T \\
B_{11}^TQ_{11}-C_{11}  & -B_{21}^TC_{22}^T-C_{13} & -C_{14} & -S_{11} & -S_{12} & -S_{13}  \\
B_{12}^TQ_{11}+Q_{31}-C_{21}  & Q_{33}-C_{23} & Q_{34}-C_{24} &  -S_{12}^T & -S_{22} & -S_{23} \\
B_{13}^TQ_{11}+Q_{41}-C_{31} & Q_{43}-C_{33} & Q_{44}-C_{34} & -S_{13}^T & -S_{23}^T & -S_{33} \rix
 \leq 0.  \label{Q}
\end{equation}
}
Let
\[  Q=X^T \mat{cccc} Q_{11} & 0 & 0  & 0\\
   0 & 0 & Q_{23} & 0         \\
  Q_{31}   & -Q_{23}^T & Q_{33} & Q_{34}\\
  Q_{41} & 0 & Q_{43} & Q_{44} \rix Y^{-1},
  \]
where by (\ref{N}) we have $ Q_{23}=-C_{22}\geq 0$.
A simple calculation, using (\ref{Q}) and $Q_{11}=\cQ_{11}\geq 0$, yields that
\[
\mat{cc} A^TQ+Q^TA & Q^TB-C^T \\ B^TQ-C & -D-D^T \rix \leq 0
\]
and $E^TQ\geq 0$.
Hence, by Theorem \ref{theorem-3}, the system (\ref{1.1}) is pH.
 \eproof

The following corollary is a consequence of Theorem \ref{theorem-4} and its proof.
\begin{corollary}Let the descriptor system (\ref{1.1}) be completely controllable, completely observable, and satisfy $D+D^T\geq 0$.
Then the system (\ref{1.1}) is positive real if and only if there exists $Q\in \Rnn$ such that   (\ref{inequality}) holds
and $E^TQ \geq 0$.
\end{corollary}

Theorem \ref{theorem-4} shows that if $D+D^T\geq 0$, then a completely controllable, completely observable and positive real system of the form (\ref{1.1}) is a pH descriptor system. But, if $D+D^T\not \geq 0$,  then the positive realness of a completely controllable and completely observable system (\ref{1.1}) does not imply that the system   is  pH, as the following example illustrates.
\begin{example}\label{ex:S} {\rm Consider the completely controllable and completely observable system  of the form (\ref{1.1}), with
\[
E=\mat{cc} 1 & 0 \\ 0 & 0 \rix, \ A=-\mat{cc} 1 & 0 \\ 0 & 1 \rix, \
   B=\mat{c} 1 \\ 1 \rix, \ C=\mat{cc} 1 & 1 \rix, \ D=-1.
\]
Then, using
\[
Q=\mat{cc} 1 & 0 \\ 0 & 1 \rix, \ W=\mat{c} 0 \\ -2 \rix
\]
one has
\[
\mat{cc} A^TQ+Q^TA & A^TW+Q^TB-C^T \\ W^TA+B^TQ-C & B^TW+W^TB -D-D^T \rix=\mat{ccc} -2 & 0 & 0 \\ 0 & -2 & 2 \\ 0 & 2 & -2 \rix \leq 0,
\]
as well as
\[
E^TQ=\mat{cc} 1 & 0 \\ 0 & 0 \rix\geq 0, \quad E^TW=0.
\]
By Theorem \ref{theorem-1}  the system (\ref{1.1})  is positive real. But, $ D+D^T=-2\not \geq 0$ and thus,  the system is not pH.
}
\end{example}
In this section we have presented a necessary and sufficient condition so that a completely controllable, completely observable, and positive real system is port-Hamiltonian. In the next section we study the case when this necessary and sufficient condition does not hold.
\section{Port-Hamiltonian realizations of positive real descriptor systems with $D+D^T\not \geq 0$}
\label{S4}
In the context of optimal control of descriptor systems different scenarios have been analyzed in \cite{ReiV19} under which solutions to Kalman-Yakubovich-Popov matrix inequalities like \eqref{1.4} and associated Lur'e equations for descriptor systems are solvable also in the case that $D+D^T\not \geq 0$. In this section we study an analogous question and derive port-Hamiltonian realizations in the case that $D+D^T\not \geq 0$. It is clear that to achieve such a realization one has to modify the feedthrough term to have a positive semidefinite symmetric part. This is only possible in the case of descriptor systems, because for standard state-space systems $E=I$, the feedthrough term is invariant for all possible realizations of the transfer function \cite{Kai80,MayA07}.

It follows from  Theorem~\ref{theorem-1} that for a completely controllable, completely observable, and positive real descriptor system (\ref{1.1})   there exists a matrix $W\in \Rnm$ such that $E^TW=0$ and $D+D^T-B^TW-W^TB\geq 0$.  The following lemma shows that there exists
such a $W$, satisfying additionally that $(E, A, C-W^TA)$ is completely observable.

\begin{lemma}\label{Lemma-W} Consider an LTI descriptor system of the form (\ref{1.1}) that is completely controllable, completely observable, and positive real. Then there exists a matrix $W\in \Rnm$ such that
\[
E^TW=0, \quad D+D^T-B^TW-W^TB\geq 0,
\]
and
$(E, A, C-W^TA)$ is completely observable.
\end{lemma}

\proof We present a constructive proof that can be implemented as a numerically reliable method.
Construct,  as in the staircase form of \cite{BunBMN99}, real orthogonal matrices $U$, $V$ and $L$ such that
\begin{eqnarray}
UEV&=&\bmat{ & \mu_1 & \mu_2 & \mu_3 \cr
\mu_1 & E_{11} & 0 & 0            \cr
\mu_2  & 0                  & 0 & 0            \cr
\mu_3  & 0                  & 0 & 0            \cr}, \quad
UAV=\bmat{ & \mu_1 & \mu_2 & \mu_3 \cr
\mu_1    & A_{11} & A_{12} & A_{13} \cr
\mu_2    & A_{21} & 0                  & 0                  \cr
\mu_3    & A_{31} & 0                  & A_{33}         \cr}, \nonumber \\
UBL&=&\bmat{ & m-\mu_2-\mu_3 & \mu_2 & \mu_3 \cr
\mu_1 & B_{11} & B_{12} & B_{13} \cr
\mu_2 & 0 & B_{22} & B_{23}                                  \cr
\mu_3 & 0 & 0 &  B_{33} \cr}, \label{EAB}
\end{eqnarray}
where $E_{11}$ and $A_{33}$ are nonsingular.  Moreover, $B_{22}$ and $B_{33}$ are also nonsingular, since $(E, A, B)$ is completely controllable which implies that
$\rank \mat{cc} E & B \rix=n$.
Set
\begin{eqnarray}
L^TCV&=&\bmat{ &\mu_1 & \mu_2 & \mu_3 \cr
m-\mu_2-\mu_3 & C_{11} & C_{12} & C_{13} \cr
\mu_2         & C_{21} & C_{22} & C_{23}  \cr
\mu_3                  & C_{31} & C_{32} & C_{33} \cr}, \nonumber \\
L^T(D+D^T)L&=&\bmat{ & m-\mu_2-\mu_3 & \mu_2 & \mu_3 \cr
m-\mu_2-\mu_1 & S_{11} & S_{12} & S_{13} \cr
\mu_2 &   S_{12}^T & S_{22} & S_{23} \cr
\mu_3 &   S_{13}^T & S_{23}^T & S_{33} \cr},  \label{CDW}\\
    UWL&=&\bmat{ & m-\mu_2-\mu_3 & \mu_2 & \mu_3 \cr
    \mu_1 & W_{11} & W_{12} & W_{13} \cr
    \mu_2 &  W_{21} & W_{22} & W_{23} \cr
\mu_3 &  W_{31} & W_{32} & W_{33} \cr}. \nonumber
 \end{eqnarray}
Since we want to have $E^TW=0$,
we choose $W_{11}=0$, $W_{12}=0$, and $W_{13}=0$.
Then
\[
L^T(C-W^TA)V=\mat{ccc} C_{11}-W_{21}^TA_{21}-W_{31}^TA_{31} & C_{12} & C_{13}-W_{31}^TA_{33} \\
C_{21}-W_{22}^TA_{21}-W_{32}^TA_{31} & C_{22} & C_{23}-W_{32}^TA_{33} \\                        C_{31}-W_{23}^TA_{21}-W_{32}^TA_{31} & C_{32} & C_{33}-W_{33}^TA_{33} \rix,
\]
and
{\scriptsize
\begin{eqnarray*}
&&L^T(D+D^T-B^TW-W^TB)L\\
&&=\mat{ccc} S_{11} & S_{12}-W_{21}^TB_{22}  & S_{13}-W_{31}^T B_{33}-W_{21}^T B_{23}  \\
S_{12}^T- B_{22}^TW_{21} & S_{22}-W_{22}^T B_{22}- B_{22}^TW_{22} &
 S_{23}-W_{32}^T B_{33}- B_{22}^TW_{23}-W_{22}^T B_{23} \\
S_{13}^T- B_{33}^TW_{31}- B_{23}^TW_{21} & S_{23}^T- B_{33}^TW_{32}-W_{23}^T B_{22}- B_{23}^TW_{22} &
S_{33}-W_{33}^T B_{33}- B_{33}^TW_{33} - B_{23}^TW_{23}-W_{23}^T B_{23}\rix.
\end{eqnarray*}
}
The existence of a matrix $W$ satisfying $E^TW=0$ and $D+D^T-B^TW-W^TB\geq 0$ implies that $S_{11}\geq 0$. Take
$ W_{23}=0$, and solve the two linear systems
\begin{equation}\label{B-1}
\mat{cc}  B_{22} &  B_{23} \\ 0 &  B_{33} \rix^T \mat{cc} W_{21} & W_{22} \\ W_{31} & W_{32} \rix
   =\mat{cc} S_{12} & S_{13} \\  {1\over 2} (S_{22}-\beta I) & S_{23} \rix^T,
\end{equation}
and
\begin{equation}\label{B-2}   B_{33}^TW_{33}={1\over 2}(S_{33}-\alpha I ),
\end{equation}
where $\beta\geq 0$,  $\alpha>0$ are chosen such that $C_{33}-W_{33}^TA_{33}$ is nonsingular, see Remark~\ref{rem:comporth} below.
We then have
\[
L^T(D+D^T-B^TW-W^TB)L=\mat{ccc} S_{11} & 0 & 0 \\ 0 & \beta I  & 0 \\ 0 & 0 & \alpha I \rix \geq 0,
\]
and hence $D+D^T-B^TW-W^TB \geq 0$.

It remains to show that $(E, A, C-W^TA)$ is completely observable. For this, set
\[
X_1=\mat{ccc} I & 0 & - A_{13} A_{33}^{-1} \\ 0 & I & 0 \\ 0 & 0 & I \rix, \quad
   Y_1=\mat{ccc} I & 0 & 0 \\ 0 & I & 0 \\
     -A_{33}^{-1} A_{31} & 0 & I \rix,
\]
so that
\begin{eqnarray*}
X_1UEVY_1&=&\mat{ccc}  E_{11} & 0 & 0 \\ 0 & 0 & 0 \\ 0 & 0 & 0 \rix, \quad
   X_1UAVY_1=\mat{ccc} \hat A_{11} & A_{12} & 0 \\  A_{21} & 0 & 0 \\ 0 & 0 &  A_{33} \rix, \\
 X_1UBL&=&\mat{ccc} \hat B_{11} & \hat B_{12} & \hat B_{13} \\ 0 &  B_{22} & B_{23} \\ 0 & 0 &  B_{33} \rix, \quad
  L^TCVY_1=\mat{ccc} \hat C_{11} &  C_{12} &  C_{13} \\ \hat C_{21} &  C_{22} &  C_{23} \\ \hat C_{31} &  C_{32} &  C_{33} \rix.
\end{eqnarray*}
Denote a full column rank matrix with its columns spanning the right  nullspace of  $E$
by $S_\infty(E)$ and with its columns spanning the right nullspace of $E^T$ by  $T_\infty(E)$, respectively. Using the fact, see \cite{ChuT08}, that a positive real system has an index at most two, it follows that in the WCF of $(E,A)$ the number of the $1\times 1$ nilpotent blocks
is $\rank(T_\infty^T(E) A S_\infty(E))$.  Next construct nonsingular matrices
$X_2$ and $Y_2$ and an  orthogonal matrix $\cL$  such that
\[
X_2 \mat{cc} E_{11} & 0 \\ 0 & 0 \rix Y_2=\bmat{ & n_1 & n_2 & n_2 \cr
n_1   & I     & 0     & 0     \cr
n_2   & 0     & 0    & I     \cr
n_2   & 0     & 0     & 0     \cr}, \quad
X_2 \mat{cc} \hat A_{11} & A_{12} \\  A_{21} & 0 \rix Y_2=\bmat{ & n_1 & n_2 & n_2 \cr
n_1 & A_1 & 0 & 0       \cr
n_2 & 0     & I & 0        \cr
n_2 & 0     & 0 & I        \cr},
\]
and that
\begin{eqnarray*}
\mat{cc} X_2 & \\ & I \rix X_1UEVY_1 \mat{cc} Y_2 & \\ & I \rix&=&\mat{cccc} I & 0 & 0 & 0 \\
0 & 0 & I & 0 \\
0 & 0 & 0 & 0 \\
0 & 0 & 0 & 0 \rix, \\
 \mat{cc} X_2 & \\ & I \rix X_1UAVY_1 \mat{cc} Y_2 & \\ & I \rix&=&\mat{cccc} A_1 & 0 & 0 & 0 \\
0 & I & 0 & 0 \\
0 & 0 & I & 0 \\
0 & 0 & 0 &  A_{33} \rix,  \\
 \mat{cc} X_2 & \\ & I \rix X_1UB L \mat{cc} \cL & \\ & I \rix&=&
\mat{ccc} \cB_{11} & \cB_{12} & \cB_{13} \\ \cB_{21} & \cB_{22} & \cB_{23} \\
              0 & \cB_{32} & \cB_{33} \\ 0 & 0 &  B_{33} \rix, \\
 \mat{cc} \cL & \\ & I \rix^T L^TCVY_1 \mat{cc} Y_2 & \\ & I \rix &=&\mat{cccc} \cC_{11} & \cC_{12} & \cC_{13} & \cC_{14} \\
\cC_{21} & \cC_{22} & \cC_{23} & \cC_{24} \\
\cC_{31} & \cC_{32} & \cC_{33} & C_{33} \rix.
\end{eqnarray*}
Since $\rank\mat{cc} E & B \rix=n$, it follows that  $\cB_{32}$ is nonsingular.
Partition the transfer function as
\[
\cT(s)=\cT_1(s)+sM,
\]
where $\cT_1(s)$ is a rational function with
$\cT_1(\infty)<\infty$.
Then
\[
M=-\mat{c} \cC_{12} \\ \cC_{22} \\ \cC_{32} \rix \mat{ccc} 0 & \cB_{32} & \cB_{33} \rix\geq 0,
\]
which implies that
\[
\cC_{12}=0, \quad \cC_{32}=0, \quad \cC_{22}\cB_{32}\leq 0, \quad \cC_{22}\cB_{33}=0,
\]
and furthermore, the complete observability of $(E, A, C)$  implies that $\rank \mat{c} E \\ C \rix=n$ and thus $\cC_{22}$ is nonsingular and hence $\cB_{33}=0$.  Using the transformation
\[
\mat{cc} X_2 & \\ & I \rix^{-T} X_1^{-T} UWL \mat{cc} \cL & \\ & I \rix=\mat{ccc} \cW_{11} & \cW_{12} & \cW_{13} \\ \cW_{21} & \cW_{22} & \cW_{23} \\
\cW_{31} & \cW_{32} & \cW_{33} \\
\cW_{41} & \cW_{42} & W_{33} \rix,
\]
together with $E^TW=0$ we obtain
\[
\cW_{11}=0, \quad \cW_{12}=0, \quad \cW_{13}=0, \quad \cW_{21}=0, \quad \cW_{22}=0, \quad \cW_{23}=0,
\]
and thus,
\[
\mat{cc} \cL & \\ & I \rix^T L^T(C-W^TA)VY_1 \mat{cc} Y_2 & \\ & I \rix
   =\mat{cccc} \cC_{11} & 0            & \cC_{13}-\cW_{31}^T  & \cC_{14}-\cW_{41}^T A_{33} \\
\cC_{21}  & \cC_{22} & \cC_{23}-\cW_{32}^T  & \cC_{24}-\cW_{42}^T A_{33} \\
\cC_{31} & 0             & \cC_{33}-\cW_{33}^T  & C_{33}-W_{33}^T A_{33} \rix.
\]
The nonsingularity of $\cC_{22}$ and  $C_{33}-W_{33}^T A_{33}$ and the complete observability of $(E, A, C)$ yield that
$(E, A, C-W^TA)$ is completely observable.
\eproof
\begin{remark}\label{rem:comporth} {\rm The matrix $C_{33}-W_{33}^TA_{33}$ is nonsingular if and only if
$\mat{cc} C_{33} & S_{33}-\alpha I \\ A_{33} & 2B_{33} \rix$ is nonsingular. Thus, using an SVD, we can compute an orthogonal matrix
\[
P=\bmat{ & \mu_3 & \mu_3 \cr
        \mu_3 & P_{11} & P_{12} \cr
        \mu_3 & P_{21} & P_{22} \cr}
\]
such that
\[
P \mat{c} C_{33} \\ A_{33} \rix=\mat{c} \Sigma_1 \\ 0 \rix,
\]
where $\Sigma_1$ is nonsingular. Note that $\mat{cc} C_{33} & I \\ A_{33} & 0 \rix$ is nonsingular,
since $A_{33}$ nonsingular,
and thus  $P_{21}$ is nonsingular. Hence
\[
P \mat{cc} C_{33} & S_{33}-\alpha I \\ A_{33} & 2B_{33} \rix=
\mat{cc} \Sigma_1 & \Lambda_1-\alpha P_{11}  \\
                 0         & \Lambda_2-\alpha P_{21} \rix.
                \]
Thus, it is easy to choose $\alpha>0$ such that $\Lambda_2-\alpha P_{21}$ is nonsingular, which gives that $C_{33}-W_{33}^TA_{33}$ is nonsingular.
}
\end{remark}

\begin{remark} \label{rem:numreliable}{\rm The matrix equations  (\ref{B-1}) and (\ref{B-2}) can be solved
via numerically reliable methods. Therefore, the numerical method for computing $W$ given in the proof of Lemma~\ref{Lemma-W} can be implemented in  numerically reliable way.}
\end{remark}

Using the construction of $W$ in Lemma~\ref{Lemma-W}, the following theorem characterizes the existence of a pH realization of  a completely controllable, completely observable, and positive real descriptor system of the form (\ref{1.1}) with $D+D^T\not \geq 0$.
\begin{theorem}\label{theorem-5}  Consider
an LTI descriptor system of the form (\ref{1.1}) that is completely controllable, completely observable, and positive real. Let $W\in \Rnm$ be such that
\[
E^TW=0, \quad D+D^T-B^TW-W^TB\geq 0,
\]
and  $(E, A, C-W^TA)$ is completely observable.
Then
\begin{equation}\label{phreal}
(E, A, B, C-W^TA, D-W^TB)
\end{equation}
is a pH  realization of system (\ref{1.1}).
\end{theorem}
\proof
It follows by direct calculation that
\begin{eqnarray*}
 \cT(s) &=& C(sE-A)^{-1}B+D \\
          &=& \mat{cc} C & I \rix \mat{cc} sE-A & 0 \\ 0 & I \rix^{-1} \mat{c} B \\ D \rix \\
          &=& \mat{cc} C & I \rix \mat{cc} sE-A & 0 \\ -sW^TE+W^TA & I \rix^{-1} \mat{c} B \\ D-W^TB \rix \\
          &=& \mat{cc} C & I \rix \mat{cc} sE-A & 0 \\  W^TA & I \rix^{-1} \mat{c} B \\ D-W^TB \rix \\
          &=& (C-W^TA)(sE-A)^{-1}B+D-W^TB.
\end{eqnarray*}
Hence, $(E, A, B,  C-W^TA,  D-W^TB)$ is a  realization of the system (\ref{1.1}), which is completely controllable and observable,  positive real, and
$ (D-W^TB) +(D-W^TB)^T \geq 0$. Therefore, by Theorem \ref{theorem-4}, $(E A, B, C-W^TA, D-W^TB)$ is port-Hamiltonian.  \eproof


After presenting conditions for the existence of a pH realization of a completely controllable, completely observable, and positive real descriptor system, in the remainder of this section, we present a computational approach to compute a  pH  realization of a general positive real descriptor system  of the form (\ref{1.1}).

Obviously, the quintuple in \eqref{phreal} is a pH realization of a system of the form (\ref{1.1}) which is completely controllable, completely observable, and positive  real.
As a result, the following corollary is a direct consequence of Theorems~\ref{theorem-4} and~\ref{theorem-5}.

\begin{corollary}\label{corollary-2} Consider a positive real LTI descriptor system of the form (\ref{1.1}). Let  $\cE$, $\cA$, $\cB$, and $\cC$ be determined by
 the  staircase form (\ref{staircase-form}).  Then following statements hold:
\begin{itemize}
\item If $D+D^T\geq 0$, then $(\cE, \cA, \cB, \cC, \cD)$ with $\cD=D$ is a pH  realization of the system (\ref{1.1}).
\item If $D+D^T\not \geq 0$, then let the matrix $\cW\in \Rnm$ satisfy
\begin{equation}\label{W2}  \cE^T\cW=0, \quad D+D^T-\cB^T\cW-\cW^T\cB\geq 0,
\end{equation}
such that
$(\cE \cA, \cC-\cW^T\cA)$ is completely observable.
Then  the quituple
\[
(\cE, \cA, \cB,  \cC-\cW^T\cA, D-\cW^T\cB)
\]
is a pH  realization of the system (\ref{1.1}).
\end{itemize}
\end{corollary}

The presented results lead to the following algorithmic framework for computing a port-Hamiltonian  realization of a general positive real system (\ref{1.1}).
\begin{enumerate}
\item [1.] Compute the  controllability-observability staircase form (\ref{staircase-form});
\item [2.] Check if $D+D^T\geq 0$.
\begin{itemize}
\item If $D+D^T\geq 0$, then $(\cE, \cA, \cB, C, \cD)$ with $\cD=D$ is a pH   realization of the system (\ref{1.1});
\item If  $D+D^T\not \geq 0$, then compute a matrix $\cW\in \Rnm$ such that (\ref{W2}) holds  and
$(\cE, \cA, \cC-\cW^T\cA)$ is completely observable.
Then with
\[      \cC:= \cC-\cW^T\cA, \ \cD= D-\cW^T\cB,
\]
the quintuple
\[ (\cE, \cA, \cB, \cC, \cD) \]
is a pH  realization of the system (\ref{1.1}).
\end{itemize}
\end{enumerate}

The  computation of the  staircase form (\ref{staircase-form}) is numerically stable provided that a proper procedure for numerical rank decisions is used, see \cite{BunBMN99}. The matrix $\cW$ in \eqref{W2} can be computed as in Lemma~\ref{Lemma-W} via
numerically reliable methods. Hence, the algorithmic framework above can be implemented in  a numerically reliable manner.

Since a completely controllable and completely observable realization of a general descriptor system can be obtained using the  staircase form (\ref{staircase-form}), in the following examples, we only consider completely controllable and completely observable descriptor systems. These examples illustrate
the  algorithmic framework.
\begin{example} \label{ex:2x2} {\rm Consider the system of the form (\ref{1.1}) in Example \ref{ex:S}, i.e.,
\[
E=\mat{cc} 1 & 0  \\ 0 & 0  \rix, \
   A=\mat{ccc} -1 & 0  \\ 0 & -1  \rix, \
 B=\mat{cc} 1 \\ 1 \rix, \ C=\mat{ccc} 1 &  1 \rix, \  D=-1.
 \]
The system is  completely controllable, completely observable, and  positive real, but it is not port-Hamiltonian.  The coefficients  $E$, $A$ and $B$ are in the form (\ref{EAB}) with $\mu_1=\mu_3=1$, $\mu_2=0$, $U=V=L=I_2$, and
\[
L^T(D+D^T)L=D+D^T=S_{33}=-2,  \quad
  UWL=W=\bmat{ &  1 \cr
1  &  W_{13} \cr
1  &  W_{33} \cr}.
\]
Following the proof of Lemma \ref{Lemma-W}, let $W_{13}=0$ and solve the matrix equation
\[
B_{33}^TW_{33}={1\over 2} (S_{33}-\alpha I), \ {\rm i.e., } \   W_{33}={1\over 2}(-2-\alpha),
\]
where $\alpha\geq 0$ and $C_{33}-W_{33}^TA_{33}=1+W_{33}^T$ is nonsingular. Take e.g. $\alpha=2$, then $W_{33}=-2$ and
\[ W=\mat{c} 0 \\ -2 \rix.
\]
With the modified
$C:=C-W^TA=\mat{cc} 1 & -1 \rix$ and $D:=D-W^TB=1$, the quintuple
$ (E, A, B, C, D) $
is a pH  realization of the system (\ref{1.1}).
}\end{example}

\begin{example} \label{ex:3x3}{\rm The LTI descriptor system of the form (\ref{1.1})  with
\begin{eqnarray*}
E&=&\mat{ccc} 1 & 0 & 0 \\ 0 & 0 & 0 \\ 0 & 0 & 0 \rix, \ A=\mat{ccc} -1 & -1 & 0 \\ -1 & 0 & 0 \\ 0 & 0 & -1 \rix, \\
  B&=&\mat{cc} 1 & 1 \\ -1 & 0 \\ 0 & 1 \rix, \ C=\mat{ccc} 1 & 1 & -1 \\ 0 & 0 & 1 \rix, \
   D=\mat{cc} -1 & 0 \\ 0 & -1 \rix.
\end{eqnarray*}
is  completely controllable, completely observable and  positive real, but it is not port-Hamiltonian. The matrices
$E$, $A$ and $B$ are in the form (\ref{EAB}) with $\mu_1=\mu_2=\mu_3=1$, $U=V=I_3$, $L=I_2$, and
\[
L^T(D+D^T)L=\bmat{ & 1 & 1 \cr
1 & S_{22} & S_{23} \cr
1 & S_{23}^T & S_{33} \cr}=\mat{cc} -2 & 0 \\ 0 & -2 \rix, \
UWL=W=\bmat{ & 1 & 1 \cr
1  & W_{12} & W_{13} \cr
1  & W_{22} & W_{23} \cr
1  & W_{32} & W_{33} \cr}.
\]
Following the proof of Lemma \ref{Lemma-W}, let $ W_{12}=0$, $W_{13}=0$, and $W_{23}=0$.
Solve the systems of linear matrix equations
\[
\mat{cc} B_{22} & B_{23} \\ 0 & B_{33} \rix\mat{c} W_{22} \\ W_{32} \rix =\mat{c} {1\over 2}(S_{22}-\beta I) \\ S_{23}^T \rix, \
   B_{33}^TW_{33}={1\over 2}(S_{33}-\alpha I),
\]
where $\beta\geq 0$, $\alpha\geq 0$ and $C_{33}-W_{33}^TA_{33}$ is nonsingular.
Inserting the coefficients we obtain the system
\[
\mat{cc} -1 & 0 \\ 0 & 1 \rix \mat{c} W_{22} \\ W_{32} \rix=\mat{c} -1-{1\over 2}\beta \\ 0 \rix, \
   W_{33}=-1-{1\over 2}\alpha,
\]
with $\beta\geq 0$, $\alpha\geq 0$, $1+W_{33}^T=-{1\over 2}\alpha$ nonsingular. With
$\beta=\alpha=2$, we obtain
\[ 
W=\mat{cc} 0 & 0 \\ 2 & 0 \\ 0 & -2 \rix. \]
With the modified
\[ C=:C-W^TA=\mat{ccc} 3 & 1 & -1 \\ 0 & 0 & -1 \rix, \
    D=: D-W^TB=\mat{cc} 1 & 0 \\ 0 & 1 \rix. \]
the quintuple
\[ (E, A, B, C, D) \]
is a pH  realization of system (\ref{1.1}).
}
\end{example}

\begin{remark}\label{rem:alternative}
{\rm A pH realization of a positive real system (\ref{1.1}) can alternatively be obtained as suggested by the following approach in  \cite{CheGH23}. Compute
$C(sE-A)^{-1}B+D=C_p(sE_p-A_p)^{-1}B_p+M_0+sM_1$,
where $E_p$  is invertible, $M_0+M_0^T\geq 0$ and $M_1=M_1^T\geq 0$. Then
\[
(\mat{ccc} E_p & 0 & 0 \\ 0 & M_1 & 0 \\ 0 & 0 & 0 \rix,~\mat{ccc} A_p & 0 & 0 \\ 0 & 0 & -I \\ 0 & I & 0 \rix, ~
     \mat{c} B_p \\ 0 \\ I \rix, ~\mat{ccc} C_p & 0 & I \rix, ~M_0)
\]
is a  pH  realization of the system (\ref{1.1}).

This realization, however, requires to evaluate and accurately compute the transfer function at $\infty$ to obtain the coefficients $M_0$ and $M_1$. Such an approach is frequently used in model reduction methods in computational fluid dynamics, see e.g. \cite{HeiSS08}. However, in the index 2 case when $M_1\neq 0$, small perturbations in the construction of $M_1$ may lead to a very different coefficient $M_1$. The perturbation analysis for this situation is rather difficult, see \cite{AhmAB10} and remains an open problem. 
}
\end{remark}

\section{Concluding Remarks}\label{S5}
The  relationship between port-Hamiltonian and positive real LTI descriptor systems has been studied. A necessary and sufficient condition is presented
 under  which a completely controllable, completely observable, and positive real LTI descriptor system is port-Hamiltonian.  A new numerically reliable
procedure is established for computing  a port-Hamiltonian  realization  of  a general positive real LTI descriptor system.
The framework is illustrated with small examples. Future work include the implementation of numerical software for the described procedures.

\end{document}